\documentclass[reqno,12pt] {amsart}

\newtheorem{theorem}{Theorem}[section]
\newtheorem{lemma}[theorem]{Lemma}
\newtheorem{corollary}[theorem]{Corollary}

\theoremstyle{definition}
\newtheorem{definition}[theorem]{Definition}
\newtheorem{assumption}[theorem]{Assumption}

\theoremstyle{remark}
\newtheorem{remark}[theorem]{Remark}

\newcommand{\mysection}[1]{\section{#1}
\setcounter{equation}{0}}

\newcommand{\bR}{\mathbb R}
\newcommand{\bS}{\mathbb S}
\newcommand{\bA}{\mathbb A}

\newcommand\cM{\mathcal{M}}

\renewcommand{\epsilon}{\varepsilon}

\begin{document}
\title[Hessian equations]{Hessian equations with
elementary symmetric functions}

\author[H. Dong]{Hongjie Dong}

\address[H. Dong]
{127 Vincent Hall, University of Minnesota, Minneapolis, MN 55455,
USA}

\email{hjdong@math.umn.edu}

\date{\today}



\begin{abstract}
We consider the Dirichlet problem for two types of degenerate
elliptic Hessian equations . New results about solvability of the
equations in the $C^{1,1}$ space are provided.
\end{abstract}

\maketitle

\mysection{Introduction}

This article is closely related to \cite{ITW}, \cite{Wang} and \cite{Kr95}. A second-order partial differential equation is called {\em Hessian}
equation if it is of the form
$$
F(u_{xx})=f,
$$
where $(u_{xx})$ is the Hessian matrix of $u$ and $F(w)$ only
depends on the eigenvalues of the symmetric matrix $w$.

Here we are concerned with the Dirichlet problem for two types of
degenerate Hessian equations:
\begin{equation}
                                      \label{eq0.1}
P_m(u_{xx})=\sum_{k=0}^{m-1}(l^+_k)^{m-k}(x)P_k(u_{xx}),
\end{equation}
\begin{equation}
                                      \label{eq0.2}
P_m(u_{xx})=g^{m-1},
\end{equation}
where $P_k(u_{xx})$ is the $k$th elementary symmetric polynomial of
eigenvalues of the matrix $u_{xx}$, and $l_k, g$ are in $C^{1,1}$.
The second equation is called $m$-Hessian equation, which becomes
the {\em Monge-Amp\`ere} equation when $m=d$. Equation \eqref{eq0.2}
is non-degenerate if $g>0$ and is degenerate if $g$ can vanish at
some points. Many authors have studied the Hessian equations and
especially the Monge-Amp\`ere equation. The solvability of {\em
non-degenerate} equations is proved by establishing the
$C^{2,\alpha}$ estimate of the solutions and using the method of continuity (see \cite{CKN1},
\cite{Kr83} and \cite{nonlinear}). For non-degenerate equations, we
also refer the reader to the works \cite{Ivoch}, \cite{CKN1},
\cite{CKN3} and \cite{Trudinger}. It is well-known that there exists
a unique admissible weak solution to {\em degenerate} problem (see
\cite{Trud3}). To show the existence and uniqueness of solutions in
$C^{1,1}$, it suffices to obtain a priori estimate of the second
order derivatives of smooth solutions to the approximating
non-degenerate equations (cf. Lemma \ref{lemma02.02}).

For the degenerate case, a global upper bound for the second-order
derivatives of admissible solutions to the Dirichlet problem of
Hessian equations, and more general {\em Bellman} equations, was
established by N.V. Krylov in a series of papers
\cite{Kr95}-\cite{payoff}. In \cite{Kr95} a few concrete equations
similar to \eqref{eq0.1} and \eqref{eq0.2} like
$$
P_m(u_{xx})=\sum_{k=0}^{m-1}(l^+_k)^{m-k+1}(x)P_k(u_{xx}),
$$
$$
P_m(u_{xx})= (l^+_k)^{m-k }(x)P_k(u_{xx}), \quad k<m,
$$
and, in particular, $P_m(u_{xx})=(g_{+})^{m}$ as $k=0$, are treated
as the applications of the general theory.

Later, for the case of the Monge-Amp\`ere equation
$$
P_d(u_{xx})=\det(u_{xx})=g^{d-1},
$$
the solvability in the $C^{1,1}$ space was proved in \cite{Wang} by
P. Guan, N.S. Trudinger and X-J. Wang with a different approach. The
power $d-1$ of $g$ was also shown to be optimal by an example in
\cite{wang2}. A modification of such example shows that $m-1$ is the
lowest possible power of $g$ when there exist second derivatives
estimates for solutions of the $m$-Hessian equations (see also
recent \cite{ITW}, where the authors did a very good survey of the
literature on the Hessian equations). As pointed out in \cite{Wang},
the techniques there rest on the fact that the solution is convex in
case $m=d$, which in general does not always hold true for
$m$-Hessian equations. With $g^{m-1}$ on the right-hand side, the
solvability of the general degenerate $m$-Hessian equations in the
$C^{1,1}$ space is still unknown.

The purpose of this paper is to prove the solvability of the
Dirichlet problems of type \eqref{eq0.1} with $C^{3,1}$ boundary
data, and also the Dirichlet problems for
  the degenerate $m-$Hessian
equations \eqref{eq0.2} with {\em homogeneous} boundary data. Our
results improve the corresponding results in \cite{Kr95}. Quite a
few arguments in the paper are based on or follow the results in
\cite{Kr95} and \cite{athm}. The technique we use is to reduce the
Hessian equations to the elliptic Bellman's equations and then apply
the general theorems on the Bellman  equations, which were
introduced in \cite{Kr95}. Owing to an observation that a certain
function is {\em quasiconvex} (Theorem \ref{thm1.4}) we are able to
apply this technique to show the solvability of  equation
\eqref{eq0.2}.

The article is organized as follows. Our two main theorems
(\ref{thm1} and \ref{thm3}) are given in the following section.
Theorem \ref{thm1} is proved in Section \ref{proofthm1}. We prove
some preliminary results and give the estimates of $u, u_x$ in
Section \ref{prelim} and \ref{est-of-ux}. In Section \ref{interior},
we use the maximum principle in a sub-domain and reduce the
estimation of second order derivatives to the estimation of
their
values on the boundary. After that, the boundary
second
derivatives are estimated in Section
\ref{mixed} in a standard way combining with a Hopf type lemma (Lemma \ref{lemma5.5}), and
this completes the proof of Theorem \ref{thm3}.

To
conclude the introduction, we explain some notation used in
what
follows: $\bR^d$ is a $d$-dimensional Euclidean space with a
fixed
orthonormal basis. A typical point in ${\mathbb R}^d$ is
denoted by
$x=(x^1,x^2,...,x^d)$. As usual the summation convention
over
repeated indices is enforced. For any
$l=(l^1,l^2,\dots,l^d)\in
\bR^d$ and any differentiable function $u$
on $\bR^d$, we denote
$D_lu=u_{x^i}l^i$ and
$D^2_lu=u_{x^ix^j}l^il^j$, etc.

Let $d\geq 2$, $m$ be positive
integers, $2\leq m\leq d$. We denote
by $\bS^d$ the set of all
symmetric $d\times d$ matrices, $\bR_+$
($\bR_+^0$) the set of all
nonnegative (strictly positive) real
numbers and $\cM_+$ ($\cM_+^0$)
the set of all nonnegative (strictly
positive) symmetric $d\times d$
matrices.

Various constants are denoted by $N$ and $\delta$ in
general and the
expression $N=N(\cdot)$ means that the given constant
$N$ depends
only on the contents of the
parentheses.

\section*{Acknowledgement}

The author would like to
express his sincere gratitude to his
advisor, N.V.~Krylov, for pointing out this problem and giving many
useful suggestions and comments for improvements. The author is also grateful to the anonymous referee who gave very useful comments on the first submission of the article.

\mysection{The Setting and Main Results}
                                          \label{mainresults}

Define $P_m(\lambda)=P_{m,d}(\lambda)$ as the $m$th elementary
symmetric polynomial of the variables
$\lambda=(\lambda^1,\lambda^2,...,\lambda^d)$. For any symmetric
$d\times d$ matrix $w$, define $\lambda(w)$ as a vector of
eigenvalues of $w$ with arbitrary order and define
$P_m(w)=P_m(\lambda(w))$. Let $C_m=C_{m,d}$ be the open connected
component of the set $\{w\in \bS^d\,:\,P_m(w)>0\}$ which contains
the identity matrix $I$.

\begin{definition}
                      \label{def1.1}
Suppose $D:=\{x\in {\mathbb R}^d \,|\, \psi(x)>0\}$ is a smooth
bounded domain with connected boundary. We say $D$ is $m-1$-convex
if for a large number $K$ and any point on $\partial D$ we have
\begin{equation}
                                  \label{eqm-convex}
P_{m-1,d-1}(\kappa^1,...,\kappa^{d-1})\geq 1/K,
\end{equation}
where $\kappa^1,...,\kappa^{d-1}$ are the principle curvatures of
$\partial D$ at this point evaluated with respect to the interior
normal to $\partial D$.
\end{definition}

\begin{definition}
A function $u\in C^2(D)$ is called $m$-admissible if the Hessian
matrix $(u_{xx})$ is in  $\bar C_m$ for any $x\in D$. A function
$u\in C^{1,1}(D)$ is called $m-$admissible if the Hessian matrix of
the second order Sobolev derivatives $(u_{xx})$ is in ${\bar C}_m$
for almost any $x\in D$.
\end{definition}

In this article, we always assume that $D$ is an $m-1$-convex bounded
domain of class $C^{3,1}$ with connected boundary. We are concerned
with the following two Dirichlet problems for the elliptic Hessian
equations,
\begin{align}
                      \label{eq1.0.0}
P_m(u_{xx})=\sum_{k=0}^{m-1}(l^+_k)^{m-k}&(x)P_k(u_{xx})\,\,\,(a.e.)
\,\,in\,\, D,\\
                      \label{eq1.0.1}
u&=\phi \,\,\,on\,\, \partial D
\end{align}
where $l_0,...,l_{m-1}$ are bounded real valued functions in $\bR^d$
and $\phi\in C^{3,1}(\bR^d)$, and the $m-$Hessian equation,
\begin{align}
                                      \label{eq1.0.2}
P_m(u_{xx})&=g^{m-1}\quad\,\,in\,\,\, D, \\
                      \label{eq1.0.3}
u&=\phi \quad\quad\,\,\,\, on\,\,\,\partial D,
\end{align}
where $2\leq m\leq d$, $g\in C^{1}(\bar D)$ is nonnegative, $\phi\in
C^{3,1}(\bar{D})$.

\begin{assumption}
                                          \label{ass1.3}
For some sufficiently large number $K$ and any $k=0,1,\cdots,m-1$
the functions $l_k(x)+K|x|^2$ and $g(x)+K|x|^2$ are convex on $\bar
D$, and
$$\text{diam}\{D\}+\|g\|_{C^{1}(\bar
D)}+1/\sup_{\bar D}|g|+\|\phi\|_{C^{3,1}(\bar
D)}+\|\psi\|_{C^{3,1}}+\|l_k\|_{C^0}\leq K,$$
$$
|\psi_x|\geq 1/K\quad \text{on}\,\partial D.
$$
\begin{equation}
                                      \label{eq1.0.0.2}
|\nabla g(x)|^2\leq Kg(x)\quad \text{in}\,D.
\end{equation}
\end{assumption}

\begin{remark}
Typically, inequality \eqref{eq1.0.0.2} holds if $g\in
C^{1,1}(\bR^d)$ is nonnegative and $\|g\|_{C^{1,1}}\leq K/2$. We
note here that \eqref{eq1.0.0.2} is also needed in \cite{Wang} as
well, for example, in the proof of Lemma 2.1 there. Because of that
the result of \cite{Wang} does not completely cover the result in \cite{Kr95}
about $\det(u_{xx})=(g_{+})^{d}$ since $(g_{+})^{d/(d-1)}$ may not
satisfy \eqref{eq1.0.0.2}.
\end{remark}

Here come our two main results.

\begin{theorem}
                                          \label{thm1}
Under the above assumptions, equation
(\ref{eq1.0.0})-(\ref{eq1.0.1}) has a unique solution $u\in
C^{1,1}(\bar{D})$ characterized by the additional property that
$$
P_m(u_{x^ix^j}+t\delta_{ij})>0,
$$
$$
P_m(u_{x^ix^j}+t\delta_{ij})>\sum_{k=0}^{m-1}(l^+_k)^{m-k}(x)
P_k(u_{x^ix^j}+t\delta_{ij})\,\,\,(a.e.) \,\,in\,\, D,
$$
for any $t>0$. Moreover, if $\sum_{k}l^+_k>0$ in $\bar{D}$, then
$u\in C^{2,\alpha}(\bar{D})$ for an $\alpha\in (0,1)$.
\end{theorem}

\begin{theorem}
                                              \label{thm3}
Under the above assumptions, there exists a unique $m-$admissible
solution $u\in C^{1,1}(\bar{D})$of the Dirichlet problem
\eqref{eq1.0.2}-\eqref{eq1.0.3} with homogeneous boundary condition
,i.e. $\phi\equiv 0$. Moreover, the solution $u$ satisfies
\eqref{eq1.0.2} almost everywhere in $D$ and admits an estimate
\begin{equation}
\|u\|_{C^{1,1}({\bar D})}\leq N(D,d,K).
\end{equation}
\end{theorem}

%For $\rho>0$, we denote $\Delta(\rho):=\{x\in D\,:\,dist(x,\partial
%D)< \rho\}$.
%\begin{theorem}
%                                            \label{thm2}
%Under
Assumption \ref{ass1.3} and in addition we assume more
%smoothness of
$g$ near $\partial D$, that is for some large $K_1$,
%\begin{equation}
%                                            \label{eq2.0.6}
%\|g^{1-1/m}\|_{C^{1,1}({\bar \Delta})}\leq K_1.
%\end{equation} Then
%equation (\ref{eq1.0.2})-(\ref{eq1.0.3}) has a unique $m-$admissible
%solution $u\in C^{1,1}(\bar{D})$ and
%\begin{equation}
%                                                \label{eq2.0.7}
%\|u\|_{C^{1,1}({\bar D})}\leq N(d,m,\rho, K,K_1).
%\end{equation}
%\end{theorem}

\mysection{Proof of Theorem \ref{thm1}}
                                      \label{proofthm1}

Firstly, let's restate some results, which can be found, for
instance, in \cite{Kr95}, as the following two lemmas.

\begin{lemma}
                                              \label{lemma1.1}
The set $C_m$ defined above is an open convex cone in $\bS^d$ with
vertex at the origin containing $\cM_+^0$, and for $k=0,1,...,m-1$
we have $C_m\subset C_k$. Furthermore, the functions
$(P_k/P_{k-1})(w),\,k=2,3,...,m$, are concave in $C_m$.
\end{lemma}

\begin{lemma}
                                              \label{lemma1.2}
If a function $F(w)$ is convex and homogeneous of degree $-\alpha$
in a cone $C\subset \bS^d$, then the function $l^{\alpha+1}F(w)$ is
a convex function of $(w,l)$ in $C\times \bR_+$.
\end{lemma}

The following lemma can be easily proved by direct calculation of
the Hessian matrix.

\begin{lemma}
                                              \label{lemma1.3}
For any positive number $n$, the function
$$
((x+y)^n-x^n)^{-1}
$$ is
a convex function in the open cone $\bR_+^0\times \bR_+^0$.
\end{lemma}

As a corollary of Lemma \ref{lemma1.2} and \ref{lemma1.3}, the
function
$$H(x,y,l):=l^{n+1}/((x+y)^n-x^n)
$$ is convex in
$\bR_+^0\times \bR_+^0\times \bR_+$.

The proof of Theorem \ref{thm1} relies on the following observation.

\begin{theorem}
                                              \label{thm1.4}
We denote
$$
G(w,l)=\sum_{k=0}^{m-1}(l_k^{m-k})(P_k/P_m)(w),
$$ where
$l=(l_0,l_1,...l_{m-1})$. Then for any
$(w,l),\,(\tilde{w},\tilde{l})\in C_m\times \bR_+^m$, we have
\begin{equation*}
G\big(\frac{w+\tilde{w}}{2},\frac{l+\tilde{l}}{2}\big)\leq
\max(G(w,l),G(\tilde{w},\tilde{l})).
\end{equation*}
That is,
$G(w,l)$ is quasiconvex in $C_m\times
\bR_+^m$.
\end{theorem}

\begin{proof}
Suppose for some $c\geq 0$,
the following two inequalities
hold
$$
\sum_{k=0}^{m-1}l_k^{m-k}P_k(w)\leq
cP_m(w),\,\,
\sum_{k=0}^{m-1}\tilde{l}_k^{m-k}P_k(\tilde{w})\leq
cP_m(\tilde{w}).
$$
Due to the homogeneity, to prove the theorem it suffices to prove
the
inequality
\begin{equation}
                                \label{eq2}
\sum_{k=0}^{m-1}(l_k+\tilde{l}_k)^{m-k}P_k(w+\tilde{w})\leq
cP_m(w+\tilde{w}).
\end{equation}
For
$k=0,...,m-1$, let $\alpha_k,\,\tilde{\alpha}_k$ be the
nonnegative
numbers such that
\begin{eqnarray}
            \label{eq1.3}
(l_k+\alpha_k)^{m-k}P_k(w)=\sum_{j=0}^{k}l_j^{m-j}P_j(w),\\
                    \label{eq1.4}
                    (\tilde{l}_k+\tilde{\alpha}_k)^{m-k}P_k(\tilde{w})
=\sum_{j=0}^{k}\tilde{l}_j^{m-j}P_j(\tilde{w}).
\end{eqnarray}
Then,
$$
(l_{m-1}+\alpha_{m-1})P_{m-1}(w)\leq
cP_m(w),\,(\tilde{l}_{m-1}+\tilde{\alpha}_{m-1})P_{m-1}(\tilde{w})\leq
cP_m(\tilde{w}).
$$
Owing
to the concavity of $(P_m/P_{m-1})(w)$, we
get
$$
c(P_m/P_{m-1})(w+\tilde{w})\geq
c(P_m/P_{m-1})(w)+c(P_m/P_{m-1})(\tilde{w})
$$
$$
\geq
l_{m-1}+\alpha_{m-1}+\tilde{l}_{m-1}+\tilde{\alpha}_{m-1},
$$
\begin{equation}
               \label{eq1.5}
cP_m(w+\tilde{w})\geq
(l_{m-1}+\alpha_{m-1}+\tilde{l}_{m-1}+\tilde{\alpha}_{m-1})
P_{m-1}(w+\tilde{w}).
\end{equation}

By
using \eqref{eq1.3} and \eqref{eq1.4}, for $k=1,...,m-1$
we
have
$$
((l_k+\alpha_k)^{m-k}-l_k^{m-k})P_k(w)=(l_{k-1}+\alpha_{k-1})^{m-k+1}P_{k-1}(w),
$$
$$
((\tilde{l}_k+\tilde{\alpha}_k)^{m-k}-\tilde{l}_k^{m-k})
P_k(\tilde{w})=(\tilde{l}_{k-1}+\tilde{\alpha}_{k-1})^{m-k+1}P_{k-1}(\tilde{w}).
$$
For
any $\epsilon>0$, the equalities above imply the following
$$
((l_k+\alpha_k+2\epsilon)^{m-k}-(l_k+\epsilon)^{m-k})
P_k(w)\geq
(l_{k-1}+\alpha_{k-1})^{m-k+1}P_{k-1}(w),
$$
$$
((\tilde{l}_k+\tilde{\alpha}_k+2\epsilon)^{m-k}-(\tilde{l}_k+\epsilon)^{m-k})
P_k(\tilde{w})\geq
(\tilde{l}_{k-1}+\tilde{\alpha}_{k-1})^{m-k+1}P_{k-1}(\tilde{w}).
$$
Again,
owing to the concavity of $(P_k/P_{k-1})(w)$ in $C_m$ and
the
corollary of Lemma \ref{lemma1.3}, we
obtain
\begin{align*}
&(P_k/P_{k-1})(w+\tilde{w})\geq
(P_k/P_{k-1})(w)+(P_k/P_{k-1})(\tilde{w})\\
&\geq
\frac{(l_{k-1}+\alpha_{k-1})^{m-k+1}}
{(l_k+\alpha_k+2\epsilon)^{m-k}-(l_k+\epsilon)^{m-k}}
+\frac{(\tilde{l}_{k-1}+\tilde{\alpha}_{k-1})^{m-k+1}}
{(\tilde{l}_k+\tilde{\alpha}_k+2\epsilon)^{m-k}
-(\tilde{l}_k+\epsilon)^{m-k}}\\
&\geq
\frac{(l_{k-1}+\alpha_{k-1}+\tilde{l}_{k-1}
+\tilde{\alpha}_{k-1})^{m-k+1}}
{(l_k+\alpha_k+\tilde{l}_k+\tilde{\alpha}_k+4\epsilon)^{m-k}
-(l_k+\tilde{l}_k+2\epsilon)^{m-k}}.
\end{align*}
As
a
consequence,
$$
((l_k+\alpha_k+\tilde{l}_k+\tilde{\alpha}_k+4\epsilon)^{m-k}
-(l_k+\tilde{l}_k+2\epsilon)^{m-k})P_k(w+\tilde{w})
$$ $$
\geq
(l_{k-1}+\alpha_{k-1}+\tilde{l}_{k-1}
+\tilde{\alpha}_{k-1})^{m-k+1}P_{k-1}(w+\tilde{w}).
$$
Letting
$\epsilon\downarrow 0$ and taking the limit
yield

$$((l_k+\alpha_k+\tilde{l}_k+\tilde{\alpha}_k)^{m-k}
-(l_k+\tilde{l}_k)^{m-k})P_k(w+\tilde{w})$$
\begin{equation}
                        \label{eq1.6}
\geq (l_{k-1}+\alpha_{k-1}+\tilde{l}_{k-1}
+\tilde{\alpha}_{k-1})^{m-k+1}P_{k-1}(w+\tilde{w}).
\end{equation}

Inequality
(\ref{eq2}) follows if we add (\ref{eq1.5}) and
(\ref{eq1.6})
together, and the theorem is proved.
\end{proof}

Finally, by relying on Theorem \ref{thm1.4} at one point, we
can essentially reproduce Krylov's approach in \cite{Kr95} with very
few modifications (cf. Remark 5.14, 5.16 of \cite{Kr95}). The idea
is to reduce the equation to a Bellman's equation and apply a
general existence and uniqueness result on degenerate elliptic
Bellman's equations proved in \cite{payoff} by a probabilistic
argument, or by an analytic approach in \cite{weak} and \cite{athm}.
In detail, due to Theorem \ref{thm1.4} we can easily get that
$$
\Theta(l):=\{w \in C_m: G(w,l)<1\}
$$
is convex in $l$ in the sense
that for any $w_i \in \Theta(l_i),i=1,2$ we have
$$
\frac{w_1+w_2}{2} \in \Theta\big(\frac{l_1+l_2}{2}\big).$$ Then
after one reduces the equation to a Bellman's equation, the free
term is semi-concave so that the regularity theory for Bellman's
equations is applicable. Theorem \ref{thm1} improves the
corresponding result in \cite{Kr95} also from the point of view of
the following remark.

\begin{remark}
The same conclusion as in
Theorem \ref{thm1} holds true if we
replace \eqref{eq1.0.0}
by
$$
P_m(u_{xx})=\sum_{k=0}^{m-1}f_k(l^+_k(x))P_k(u_{xx})\,\,\,(a.e.)
\,\,in\,\,
D,
$$
for any $s\geq 0$, where $f_k:\bR_+\to \bR_+$ are
continuous
functions and $f_k^{1/(m-k)}$ are convex for
$k=0,1,\cdots,m-1$.

Indeed, if we
denote
$$
G_1(w,l)=\sum_{k=0}^{m-1}f_k(l_k) (P_k/P_m)(w),
$$ then due
to Theorem \ref{thm1.4} and the convexity of
$f_k^{1/(m-k)}$ in
$\bR_+$, for any
$(w,l),\,(\tilde{w},\tilde{l})\in C_m\times \bR_+^m$
we
have
$$
\max(G_1(w,l),G_1(\tilde{w},\tilde{l}))
$$
$$
=\max\big\{\sum_{k}f_k(l_k)(P_k/P_m)(w),\sum_{k}
f_k(\tilde{l}_k)
(P_k/P_m)(\tilde{w})\big\}
$$
$$
\geq
\sum_{k}\big[\frac{1}{2}(f_k^{1/(m-k)}(l_k)
+f_k^{1/(m-k)}(\tilde{l}_k))\big]^{m-k}
(P_k/P_m)\big(\frac{w+\tilde{w}}{2}\big)
$$
$$
\geq
\sum_{k}f_k\big(\frac{l_k+\tilde{l}_k}{2}\big)
(P_k/P_m)\big(\frac{w+\tilde{w}}{2}\big)=
G_1\big(\frac{w+\tilde{w}}{2},\frac{l+\tilde{l}}{2}\big),
$$
i.e.
$G_1(w,l)$ is also quasiconvex in $C_m\times \bR_+^m$. Thus
our
assertion follows.
\end{remark}

\mysection{Some Preliminary
Results}

\label{prelim}

In what follows, we consider the following Dirichlet
problem for
$2\leq m\leq d$:
\begin{align}
                          \label{dirichlet}
P_m(u_{xx})&=g^{m-1}\quad\,\,in\,\,\, D, \\
                          \label{boundary}
u&=0 \quad\quad\,\,\,\, on\,\,\,\partial D,
\end{align}
where $g\in C^{1,1}(\bar{D})$ is nonnegative and $P_m$
is defined in
Section \ref{mainresults}. We focus on the solution $u$
such that
$(u_{xx}) \in {\bar C}_m$ for any $x\in
D$.

\begin{lemma}
                            \label{lemma2.1}
For any matrix $(v_{ij})$ in $C_m$, the $d$ by $d$ matrix
$K(v):=\big(P_{m,v_{ij}}(v)\big)$ is positive definite.
\end{lemma}
\begin{proof}
Let $\eta$ be a nonzero vector in
$\bR^d$. We have
$$
\eta^TK(v)\eta=\text{Tr}(K(v)\eta\eta^T)
=\frac{d}{dt}P_m(v+t\eta\eta^T)\big|_{t=0}.
$$
Since
the rank of $\eta\eta^T$ is one, $P_m(v+t\eta\eta^T)$ is
linear in
$t$. Thus,
$$
P_m(v+t\eta\eta^T)=P_m(v)+t \eta^TK(v)\eta.
$$
Also
note that $\eta\eta^T\in \bar{C}_d \subset \bar{C}_m$, for any
$t\geq 0$
we have
$$
v+t\eta\eta^T\in \bar{C}_m,\quad P_m(v+t\eta\eta^T)\geq
0.
$$
Therefore, $\eta^TK(v)\eta\geq 0$ is nonnegative.
If
$\eta^TK(v)\eta=0$, for any $t> 0$ it holds
that
$$
P_m(v-t\eta\eta^T)=P_m(v)>0,
$$
which implies
$v-t\eta\eta^T\in C_m$ and thus $ v/t-\eta\eta^T\in
C_m $. Letting
$t\to +\infty$ yields $-\eta\eta^T\in {\bar C}_m\in
{\bar C}_1$. But
this is impossible
because
$$
P_1(-\eta\eta^T)=\text{Tr}(-\eta\eta^T)=-|\eta|^2<0.$$

Hence
it holds that $\eta^TK(v)\eta>0$, and the lemma is
proved.
\end{proof}

The following corollary is an immediate
consequence of Lemmas
\ref{lemma1.1} and
\ref{lemma2.1}.

\begin{corollary}
                                \label{cor4.2}
For any $d\times d$ matrix $(v_{ij})$ in $C_m$, all the $(d-1)\times
(d-1)$ submatricies obtained by deleting the $k$th row and $k$th
column of $(v_{ij})$ ($k=1,\cdots,d$) are in $C_{m-1,d-1}$.
\end{corollary}

\begin{lemma}
                    \label{lemma4.2.1}
(i) For any matrices $(v_{ij})$ and $(w_{ij})$ in $C_m$, we have
$$\text{Tr}(K(v)w)=P_{m,v_{ij}}(v)w_{ij}>0.$$

(ii)
Moreover, for any orthogonal matrix $Q$, it holds
that
\begin{equation}
                                    \label{eq4.1.10.57}
K(QvQ^{T})=QK(v)Q^{T}.
\end{equation}
\end{lemma}
\begin{proof}
The key idea of the proof of part (i) is to use the properties of
hyperbolic polynomials. We treat $P_m(v)$ as a homogenous polynomial
of $d(d+1)/2$ variables
$$v_{11},v_{12},...
,v_{1d},v_{22},v_{23},...,v_{2d},v_{33},...,v_{3d},...,v_{dd}.$$
It's known that for any $\mu\in \bR^d$, all roots of the polynomial
$P_m(\mu+t\lambda_0)$ are real, where $\lambda_0=(1,1,...,1)$ (cf.
Corollary 6.5 of \cite{Kr95}). Therefore, by the very definition of
$P_m$ in $\bS^d$, for any $v\in \bS^d$, all roots of the polynomial
$P_m(v+tI_d)$ are real. Owing to Lemma 4.16 and Theorem 6.4 of
\cite{Kr95}, for any $v, w$ in $C_m$ we have
$$
P_{m,v_{ij}}(v)w_{ij}=\frac{d}{dt}P_m(v+tw)\big|_{t=0}>0.
$$
This proves the first part of the lemma.

To prove the second part, first notice that $P_m$ and $C_m$ are both
invariant under orthogonal transformations. Suppose $w$ is a
symmetric matrix. We have
$$
\text{Tr}(K(QvQ^T)w)=\frac{d}{dt}P_m(QvQ^T+tw)\big|_{t=0}=
\frac{d}{dt}P_m(v+tQ^TwQ)\big|_{t=0}
$$
$$
=\text{Tr}(K(v)Q^TwQ)=\text{Tr}(QK(v)Q^Tw).
$$ Since $w$ is an arbitrary symmetric matrix, the conclusion of
(ii) follows immediately.
\end{proof}

\begin{remark}
					\label{remark11.15}
By using the same method as in the proof of \ref{lemma4.2.1} (i), owing to Lemma 4.16 and Theorem 6.4 of \cite{Kr95}, one can prove that
for any $(v_{ij}),(w_{ij})$ in $C_m$ and $t_0,s_0\geq 0$ it holds that
$$
\frac{\partial^2}{\partial t\partial s}P_m(sv+tw)\big|_{(s,t)=(s_0,t_0)}\geq 0.
$$
And this implies 
\begin{equation}
				\label{eq12.10}
P_m(v+w)\geq P_m(v)+P_m(w).
\end{equation}
\end{remark}
\begin{lemma}
                                     \label{lemma02.04.2}
Let $v\in \bar{C}_m$ and $c$ be a nonnegative constant. Then
$P_m(v)=c$ if and only if
\begin{equation}
                                     \label{eq02.04.4}
\inf_{w\in
C_m}\{a^{ij}(w)v^{ij}-m(d-m+1)^{-1}P_m^{1-1/m}(w)
P_{m-1}^{-1}(w)c^{1/m}\}=0,
\end{equation}
where
\begin{equation}
                                              \label{eq2.3}
a^{ij}(w)=P_{m,w^{ij}}(w)/\text{Tr}(K(w)).
\end{equation}
\end{lemma}

\begin{proof}
First we suppose $P_m(v)=c$. Recall that $P_m^{1/m}(v)$ is concave
in $C_m$ (see, for instance,  Theorem 6.4 of \cite{Kr95}). We get
$$
c^{1/m}=P^{1/m}_m(v)
$$
$$
=\inf_{w\in C_m}\{m^{-1}P_{m,w^{ij}}(w)P_m^{1/m-1}(w)(v^{ij}-w^{ij})
+P_m^{1/m}(w)\}
$$
\begin{equation}
                                          \label{eq2.2}
=\inf_{w\in
C_m}\{m^{-1}P_{m,w^{ij}}(w)P_m^{1/m-1}(w)v^{ij}(x)\}.
\end{equation}
The last equality is because $P_m$ is a homogeneous polynomial of
degree $m$. In case $c>0$, because the infimum above is attained
when $\omega=v$ and also because $\text{Tr}(K(v))>0$, equality
\eqref{eq02.04.4} follows immediately. If $c=0$, then one has $v\in
\partial C_m$. Note that
$$
\text{Tr}(K(v))=(d-m+1)P_{m-1}(v)
$$
and $v+\epsilon I\in C_m$. For any real number $\epsilon>0$,
we have
$$
a^{ij}(v+\epsilon I)v^{ij}=P_{m,v^{ij}}(v+\epsilon I)v^{ij}/
\text{Tr}(K(v+\epsilon I))
$$
$$
=P_{m,v^{ij}}(v+\epsilon I)(v^{ij}+\epsilon \delta_{ij}-\epsilon
\delta_{ij})/ \text{Tr}(K(v+\epsilon I))
$$
$$
=m(d-m+1)^{-1}(P_m/P_{m-1})(v+\epsilon I)-\epsilon.
$$
Thanks to Lemma 4.16 (ii) of \cite{Kr95}, we get
$$
\lim_{\epsilon\downarrow 0}a^{ij}(v+\epsilon I)v^{ij}=0,
$$
and \eqref{eq02.04.4} follows.

On the other hand, assume that
  \eqref{eq02.04.4} holds true. Observe that
$$
P_m(v+tI)\to +\infty\quad \text{as}\,\, t\to +\infty
$$
and there exists $t_0\leq 0$ such that $v+t_0I\in \partial C_m$ and
$P_m(v+t_0I)=0$. By continuity one can find a real number
$t_1\geq
t_0$ such that $v+t_1I\in {\bar C}_m$ and $P_m(v+t_1I)=c$.
Due to
the first part of the proof, it holds that
$$
\inf_{w\in
C_m}\{\text{Tr}(a(w)(v+t_1I))-m(d-m+1)^{-1}P_m^{1-1/m}(w)
P_{m-1}^{-1}(w)c^{1/m}\}=0,
$$
$$
\inf_{w\in C_m}\{\text{Tr}(a(w)v)-m(d-m+1)^{-1}P_m^{1-1/m}(w)
P_{m-1}^{-1}(w)c^{1/m}\}=-t_1.
$$
Therefore, from \eqref{eq02.04.4} we obtain $t_1=0$, and the lemma
is proved.
\end{proof}

Because of Lemma \ref{lemma02.04.2}, equation \eqref{dirichlet} is
equivalent to
\begin{equation}
                                          \label{eq2.2.b}
\inf_{w\in C_m}\{L^{w}u(x)+f(w,x)\}=0,
\end{equation}
where
$$
L^{w}u(x)=a^{ij}(w)u_{x_ix_j}(x),
$$
$$
f(w,x)=-m(d-m+1)^{-1}P_m^{1-1/m}(w)P_{m-1}^{-1}(w)g^{1-1/m}(x).
$$

Owing to Theorem 1.1 of \cite{Trud3},
\eqref{dirichlet}-\eqref{boundary} has a unique admissible weak
solution $u\in C^0(D)$ in the sense that for a sequence of
$m$-admissible functions $u_k\in C^2(D)$ we have
$$u_k\to u \,\,\text{in}\,\,C^0(D),\quad P_m(u_{k,xx})
\to g^{m-1}\,\,\text{in}\,\,L_{\text{loc}}^1(D).$$

The following lemma will be proved in Section \ref{prooflemma}.
\begin{lemma}
                     \label{lemma02.02}
Under Assumption \ref{ass1.3}, if we can establish a priori $C^2$
estimate of solutions to non-degenerate problems with $C^2(\bar D)\in
g>0$, which does not depend on the infimum of $g$ in $D$, then

(i) the admissible weak solution $u$ of
\eqref{dirichlet}-\eqref{boundary} is in $C^{1,1}(\bar D)$.

(ii) Moreover, $u$ satisfies \eqref{dirichlet} almost everywhere in
$D$.
\end{lemma}

Due to Lemma \ref{lemma02.02}, from now on, we always assume that
$g$ is positive on $\bar D$ and belongs to $C^2(\bar D)$.
  In this case, it's known that $u\in C^4(D)\cap
C^2({\bar D})$. Because the infimum in \eqref{eq2.2.b} is attained
when $\omega=u_{xx}$, it is easy to see that \eqref{eq2.2.b} is
equivalent to a uniformly elliptic Bellman equation
\begin{equation}
                                          \label{eq05.02.11}
\inf_{w\in C_m^*}\{L^{w}u(x)+f(w,x)\}=0,
\end{equation}
where $$C_m^*=\{w\in C_m|P_m(w)\geq \inf_Dg^{m-1},w\leq
N(\|u\|_{C^2({\bar D})})I\}$$

Let $\xi$ be a unit vector in $\bR^d$. Denote
$\lambda^1(x),...,\lambda^d(x)$ to be the eigenvalues of the Hessian
matrix $u_{xx}(x)$. Due to Lemma \ref{lemma1.1}, we have
\begin{equation}
                                      \label{traceest}
\lambda^1+\lambda^2+...+\lambda^d=\text{Tr}(u_{xx})=\Delta
u=P_1(u_{xx})> 0.
\end{equation}

We can get more than (\ref{traceest}). Define $\Lambda_2$ as the
open connected component of the set $\{\lambda\in
\bR^d:P_2(\lambda)>0\}$ which contains the vector $(1,1,...,1)$.
Since $(u_{xx})\in C_m\subset C_2$, we have
$\lambda:=(\lambda^1,\lambda^2,...,\lambda^d)\in \Lambda_2$.
Obviously, $\lambda_{\epsilon}:=(1,\epsilon,\epsilon,...,\epsilon)$
is also in $\Lambda_2$ for any $\epsilon>0$. Thus by Theorem 6.4(i)
in \cite{Kr95} with $\lambda_{\epsilon}$ and $P_2$ in place of
$\lambda_1$ and $Q_m$ respectively, we get
$$
0<\frac{\partial}{\partial
t}P_2(\lambda+t\lambda_{\epsilon})=\lambda^2+\lambda^3+...+\lambda^d+O(\epsilon).
$$
By letting $\epsilon\downarrow 0$ and taking the limit, we get
\begin{equation}
                                          \label{eq2.05}
\lambda^2+\lambda^3+...+\lambda^d=\text{Tr}(u_{xx})-\lambda^1\geq
0.
\end{equation}
Of course, (\ref{eq2.05}) remains true if we replace $1$ by any
$i=2,3,...,d$.

If $\max_i\lambda^i(x)<1$, due to (\ref{traceest}), we immediately
get an estimate of $\lambda^i(x),\,i=1,2,...,d$. And this yields the
estimate of $(u_{xx})$. Therefore, in the sequel we only consider
the region
$$
D':=\{x\in D|\,\max_i\lambda^i(x)\geq1\}.
$$
Due to (\ref{eq2.05}), in $D'$ we have $P_1(u_{xx})\geq 1$.

\begin{lemma}
                                  \label{lemma01.03}
(i) For any $x\in D'$, we have
\begin{equation}
                                          \label{eq2.7}
g^{m-2}(x)\leq \text{Tr}(K(u_{xx}))/(d-m+1).
\end{equation}
(ii) For any $x\in D$, we have
\begin{equation}
                                          \label{eq2.7.1}
g^{m-2}(x)\leq Ng^{m-3/2}(x)\leq N\text{Tr}(K(u_{xx}))/(d-m+1),
\end{equation}
where $N$ depends only on $K$.
\end{lemma}

\begin{proof} Note that $(P_k(w)/{d \choose k})$ is a
log-concave function of $k$ for $k=0,1,...,m$ (cf., for instance,
  Corollary 6.5  \cite{Kr95}). So in $D'$ we have
$$
g^{m-2}(x)=P_m^{(m-2)/(m-1)}(u_{xx})\leq
P_m^{(m-2)/(m-1)}(u_{xx})P_1^{1/(m-1)}(u_{xx})
$$
$$
\leq NP_{m-1}(u_{xx})=N\text{Tr}(K(u_{xx}))/(d-m+1).
$$
Also for any
matrix $w\in C_m$,
\begin{equation}
                                          \label{eq2.7.2}
P_m^{1-1/m}(w)=P_0^{1/m}(w)P_m^{(m-1)/m}(w) \leq NP_{m-1}(w).
\end{equation}
Thus in $D$, we have
$$
g^{m-3/2}(x) \leq Ng^{(m-1)^2/m}(x)= NP_m^{(m-1)/m}(u_{xx}(x))
$$
$$
\leq NP_{m-1}(u_{xx}(x)) =N\text{Tr}(K(u_{xx}))/(d-m+1).
$$
The lemma is proved
\end{proof}

Due to Assumption \ref{ass1.3} and \eqref{eq2.7.2}, for any $w\in
C_m$ we have
\begin{equation}
                      \label{eq2.0.0.4}
\|f(w,\cdot)\|_{C^1({\bar D})}\leq N(K,d).
\end{equation}

Near $\partial D$ let $\psi_2=\text{dist}(x,\partial D)$ if $x\in
\bar{D}$, $\psi_2=-\text{dist}(x,\partial D)$ if $x\notin D$. By a
standard argument in \cite{Kr95}, we can define
$\psi_1=\psi_2-t\psi_2^2$ near $\partial D$ with $t$ sufficiently
large and continue $\psi_1$ in an appropriate manner such that
$$
\{x\in \bR^d\,:\,\psi_1(x)>0\}=D,\,\,
\|\psi_1\|_{C^{3,1}(\bR^d)}\leq N(K,d)
$$ and we can find
$\rho=\rho(K,d)>0$ such that for any $x\in \bar{\Delta}_{\rho}$,
$-\psi_{1,xx}$ is in $C_m$, where
$$\Delta_{\rho}=\{x\in D\,:\,\text{dist}(x,\partial D)<\rho\}.$$

Denote $\Omega$ to be the  closure of
$\{K(v)/\text{Tr}(K(v))\,|\,v\in C_m\}$. Obviously, $\Omega$ is a
compact set. Owing to Lemma \ref{lemma4.2.1}, by a compactness
argument we obtain the following corollary.
\begin{corollary}
                                          \label{cor2.2}
There exist a $\delta=\delta(K,d,\rho)>0$, such that for any $v
\in
\Omega$ and $x\in \bar{\Delta}_{\rho}$ we
have
$\text{Tr}(v\psi_{1,xx}(x))<-\delta$. Especially, for any $w\in
C_m$
and $x\in \bar{\Delta}_{\rho}$, we have
$L^{w}\psi_1(x)<-\delta$.
Moreover, for any $x\in \partial D$,
$|\psi_{1x}|=1$.
\end{corollary}

\begin{proof}
For any $x\in
\bar{\Delta}_{\rho}$ we can find a number $\epsilon>0$
such that
$-\psi_{1,xx}-\epsilon I\in C_m$. Then we
have
$$
\text{Tr}(v(-\psi_{1,xx}-\epsilon I))\geq
0,\quad
\text{Tr}(v(-\psi_{1,xx}))=\epsilon> 0.
$$
Because both
$\Omega$ and $\bar{\Delta}_{\rho}$ are compact, there
exists
$\delta>0$ such that $\text{Tr}(v\psi_{1,xx}(x))<-\delta$.
The last
assertion follows immediately from the definition
of
$\psi_1$.
\end{proof}

Let $\psi_0=(2R+1)^2-|x-x_0|^2$, where $R$
is the diameter of $D$
and $x_0$ is a point in $D$. Since
$\text{Tr}(a(w))=1$, it holds
that $L^w\psi_0=-2$ for any $w\in C_m$.
Thus Assumption 1.2 (a), (b)
in \cite{payoff} are satisfied and the
following lemma is proved in
\cite{payoff}, Lemma 1.1.
\begin{lemma}
                       \label{lemma2.3.1}
There exist $\psi\in C^{3,1}$ and $\delta=\delta(K,d)>0$ satisfying
(i) $\{x\in \bR^d\,:\, \psi(x)>0\}=D$, (ii) $|\psi_x|>\delta$ on
$\partial D$, (iii) $\|\psi\|_{C^{3,1}(\bR^d)}\leq N(K,d)$ and (iv)
for any $w \in C_m$ and $x\in D$ we have $L^w\psi(x)<-\delta$.
\end{lemma}

\mysection{Estimates of $u$ and
$u_x$}
                                     \label{est-of-ux}

In
this section, we will give some estimates of $u$ and
$u_x$.
\begin{theorem}
                         \label{theorem3.1}
For
any $x\in D$, we have
$$|u(x)-\phi(x)|\leq
N(K,d)\psi(x).$$
\end{theorem}
\begin{proof}
Denote
$\tilde{u}=u-\phi$. Since $P_m$ is a homogeneous polynomial
of degree
$m$, we
have
$$
P_{m,u_{x^ix^j}}(u_{xx}(x))u_{x^ix^j}(x)=mP_m(u_{xx}(x))=mg^{m-1}(x),
$$
\begin{equation}
                              \label{eq3.1}
\text{Tr}\big(a(u_{xx}(x))\tilde{u}_{xx}\big)
=mg^{m-1}(x)/\text{Tr}\big(K(u_{xx}(x))\big)-
\text{Tr}\big(a(u_{xx}(x))\phi_{xx}\big).
\end{equation}

 From the positiveness of $K(u_{xx})$, equality (\ref{eq3.1}) can
be
looked at as a second order elliptic equation of $\tilde{u}$. Due
to
(\ref{eq2.7.1}), the right hand side of (\ref{eq3.1}) is bounded
by
a constant depending only on $K$, $d$ and $m$.
Lemma
\ref{lemma2.3.1}
implies
$$
\text{Tr}\big(a(u_{xx}(x))\psi_{xx}(x)\big)<
-\delta.
$$
After using the comparison principle, we get what we
expected.

\end{proof}

\begin{theorem}
                     \label{theorem3.2}
We have $\sup_D|u_x|\leq N$, where $N$ depends only on $K$ and $d$.
\end{theorem}
\begin{proof}
After
differentiating \eqref{dirichlet} in the direction $\xi$,
we
get
\begin{equation}
                                \label{eq3.3.3}
a^{ij}(u_{xx}(x))u_{(\xi)x^ix^j}(x)
=\frac{m-1}{d-m+1}g_{(\xi)}(x)g^{m-2}(x)P_{m-1}^{-1}(u_{xx}(x)).
\end{equation}
Owing
to \eqref{eq1.0.0.2} and \eqref{eq2.7.1}, the absolute value
of the
right-hand side of \eqref{eq3.3.3} is less
than
$$
N|g_{(\xi)}(x)|g^{m-2}(x)P^{-1}_{m-1}(u_{xx}(x))
\leq
Ng^{m-3/2}P^{-1}_{m-1}(u_{xx}(x))|\leq N.
$$
Again by using
Lemma \ref{lemma2.3.1} and the comparison principle,
we
get
\begin{equation*}
\sup_D|u_{(\xi)}|\leq N+\sup_{\partial
D}|u_{(\xi)}|.
\end{equation*}

Upon using Lemma \ref{theorem3.1}, we
get the estimate of the first
derivative on the
boundary:
$$|u_{(\xi)}(x)|\leq N(K,d)\quad \forall
x
\,\,\text{on}\,\, \partial D,$$ and the lemma is
proved.
\end{proof}

We also need a lower bound for the normal first
order derivative
$D_nu$ on the boundary $\partial D$.
\begin{lemma}
                          \label{lemma5.5}
Under Assumption \ref{ass1.3}, for any $x\in
\partial D$ we have
$$
D_nu(x)\leq -\gamma<0,
$$ where $n$ is the
inner normal of $\partial D$ at point $x$ and
$\gamma$ is a constant
depending only on $D$, $K$ and $d$.
\end{lemma}
\begin{proof}
 From
\eqref{eq2.2} with $w=I_d$ we get
\begin{equation}
                                    \label{eq5.5.5}
\Delta u(x)\geq N(m,d)g^{1-1/m}(x).
\end{equation}
Because $g$ is not
always equal to $0$ in $D$ and $u=0$ on $\partial
D$, by the strong
maximum principle, we get $u$ is strictly positive
in $D$.

Since $D$ is a $C^{3,1}$ domain, we can find two positive
numbers
$\epsilon_1$ and $\epsilon_2$ depending only on $K$
satisfying: (i)
for any $x\in
\partial D$, there is a ball $B_x$ of
radius $\epsilon_1$ which is inside
$D$ and $\partial B_x$ and
$\partial D$ are tangent at $x$; (ii) for
any $x\in \partial D$, if
we look at $x$ as "north poll", then the
south half ball of $B_x$ is
in $D\setminus \Delta_{\epsilon_2}$.

Next, we claim that there
exists $\gamma_0=\gamma_0(D,d,K)>0$ such
that
$$
u(x)\leq
-\gamma_0\quad\, \text{in}\,D\setminus \Delta_{\epsilon_2}.
$$
We
prove this by contradiction. If this is not true, we can find
two
sequences $x_j\in {\bar D}$, $g_j\in C^{2}({\bar D})$ such
that
\begin{equation}
                                    \label{eq5.5.7}
\|g_j\|_{C^{1}({\bar D})}\leq K,\quad g_j(x_j)\geq 1/K,
\end{equation} and
also $u_j$ solves
\eqref{dirichlet}-\eqref{boundary} with $g_j$ in
place of $g$ such
that
\begin{equation}
                                    \label{eq5.5.6}
\sup_{D\setminus \Delta_{\epsilon_2}}u_j\to 0,\quad \text{as}\,\,
j\to +\infty.
\end{equation} Because $\bar{D}$ is
compact, after passing to a
subsequence, we may assume $x_j$
converges to a point $x_0\in {\bar
D}$. Due to \eqref{eq5.5.7} for a
small neighborhood $U_{x_0}$ of
$x_0$ we have $\inf_{U_{x_0}\cap
{\bar D}}g_j\geq 1/(2K)$ for all
$j$ large enough. Let $u_0$ be the
solution of
$$
\Delta u_0(x)=N(m,d)(I_{U_{x_0}\cap
D}/(2K))^{1-1/m}(x)
$$
with zero boundary data, where $N(m,d)$ is the
same positive
constant as that of \eqref{eq5.5.5}. Then by the
comparison principle, we get
$$
u_j(y)\leq u_0(y)\quad
\,\text{in}\,\,D.
$$ However, by the strong maximum
principle, $u_0$
is strictly negative in $D$ and bounded away from $0$ on $D\setminus
\Delta_{\epsilon_2}$, which contradicts
\eqref{eq5.5.6}.

Now for any
$x\in \partial D$, we consider $u$ in the ball $B_x$. By
the previous
proof, we have $u\leq -\gamma_0$ on the south half
sphere and
$u(x)=0$. Moreover, $u$ is subharmonic in $B_x$. Denote
$v$ to a
harmonic function in $B_x$ with boundary data $0$ on the
north half
sphere and $-\gamma_0$ on the south half sphere. By the
comparison
principle again we have $u(y)\leq v(y)$ in $B_x$, and
therefore
$D_nu(x)\leq D_nv(x)<0$. Actually, by our construction of
$v$,
$D_nv(x)$ only depends on $K$ and $d$ (not $x$). This completes
the
proof of the lemma.
\end{proof}

\mysection{Interior second order derivatives} \label{interior}

Here, our goal is to firstly give an interior
estimate of the second
order derivatives of the solution via the
estimates on the boundary
of the second order derivatives.

Note that
for any function $H(\alpha,x)$ which is twice
differentiable in $x$,
if $\inf_{\alpha\in A}H(\alpha,x)$ is also
twice differentiable, then
for any $\xi\in \bR^d$ pointwisely we
have
$$\big(\inf_{\alpha\in
A}H(\alpha,x)\big)_{(\xi)(\xi)}\leq
H_{(\xi)(\xi)}(\alpha_0,x),\quad
\Delta\big(\inf_{\alpha\in
A}H(\alpha,x)\big)\leq \Delta
H(\alpha_0,x),$$ where $\alpha_0\in A$
such that
$H(\alpha_0,x)=\inf_{\alpha\in A}H(\alpha,x)$. After
differentiating
(\ref{eq2.2.b}) twice in the direction $\xi$, we
get
$$
(P_m^{1-1/m}P^{-1}_{m-1})(u_{xx}(x))\frac{m-1}{d-m+1}
\big[g_{(\xi)(\xi)}g^{-1/m}(x)
$$
$$-(1/m)g_{(\xi)}^2g^{-1-1/m}(x)\big]\leq
a^{ij}(u_{xx}(x))u_{x^ix^j(\xi)(\xi)}(x).
$$
Because
of \eqref{eq1.0.0.2}, \eqref{dirichlet} and \eqref{eq2.3},
we
obtain
\begin{equation}
                                    \label{eq2.4}
P_{m,u_{x^ix^j}}(u_{xx}(x))u_{x^ix^j(\xi)(\xi)}(x)\geq -Ng^{m-2}(x).
\end{equation}

Observe that by Lemma \ref{lemma2.1}, (\ref{eq2.4}) is an elliptic
equation in $D$. Combining (\ref{eq2.7}) with the comparison
principle for the elliptic equations in $D'$, we get an upper
estimate
\begin{equation}
                                          \label{uppest}
u_{(\xi)(\xi)}(x)\leq N+\sup_{\partial D'}u_{(\xi)(\xi)}(x)\leq
N+\sup_{\partial D}u_{(\xi)(\xi)}(x).
\end{equation}
The last inequality is because $|u_{(\xi)(\xi)}(x)|\leq d-1$ on
$\partial D'\cap D$.

To obtain the lower estimate, it remains to use (\ref{traceest})
again. As a conclusion, we get
\begin{theorem}
                                          \label{theorem2.1}
Let $u\in C^4(D)\cap C^2(\bar{D})$ be the solution of
(\ref{dirichlet}) in $D$ and satisfies $(u_{xx}(x))\in C_m$ in $D$.
Then
\begin{equation}
                                          \label{eq2.09}
\sup_{D}|D^2u|\leq N(1+\sup_{\partial D}|D^2u|),
\end{equation}
where N is a constant depending only on $K$ and $d$.
\end{theorem}

\begin{remark}
It turns out that to get \eqref{eq2.09} it suffices to assume
$\Delta g$ to be bounded from below. Indeed, under this condition,
instead of \eqref{eq2.4} and \eqref{uppest} we have
$$
P_{m,u_{x^ix^j}}(u_{xx}(x))\Delta u_{x^ix^j}(x)\geq -Ng^{m-2}(x),
$$
$$
\Delta u(x)\leq N+\sup_{\partial D'}\Delta u(x)\leq N+\sup_{\partial
D}\Delta u(x).
$$ The last estimate together with \eqref{eq2.05} yields
\eqref{eq2.09}.
\end{remark}

\mysection{Boundary second order derivatives}
                                          \label{mixed}

We remark here that for the problem with homogeneous boundary condition the estimate of the second order derivatives on the boundary follows from the arguments in Section 5 of \cite{CKN3}\footnote{We are grateful to the referee for
pointing this out.} by applying Lemma \ref{lemma5.5} instead of the usual Hopf lemma. Here we give some details for the sake of completeness. While estimating the mixed second derivates, we use Krylov's approach in \cite{weak} and \cite{athm}.

For any $x\in \partial D$, after a shift of the origin and an
orthogonal transformation, we may suppose $x$ is the origin and
$x^n$-axis is the inner normal. By further transforming the
coordinate $x'=(x^1,...,x^{n-1})$, we can assume in a small
neighborhood $U_0$ of $x$, $\partial D$ can be represented by
$x^n={\bar \psi}(x')$ and $u_{x_ix_j}(0)=0$ for $i\neq j$,
$i,j=1,...,n-1$. Here $x^n-{\bar \psi}(x')$ is in the class of
$C^{3,1}({\bar U}_0)$ and
$$
\|x^n-{\bar \psi}(x')\|_{C^{3,1}({\bar U}_0)}\leq N(d,K),\quad
\nabla {\bar \psi}(0)=0.
$$
Then it suffices to estimate $u_{x^jx^j}(0)$, $u_{x^jx^n}(0)$ and
$u_{x^nx^n}(0)$, where $j=1,...,n-1$.

The estimation of the tangential second order derivatives on the
boundary is standard (cf. \cite{CKN3}, \cite{Wang} or \cite{athm}). We
differentiate the equality $$u(x',{\bar \psi}(x'))=0$$ twice
with
respect to $x^j$, $j=1,...,n-1$, and get
\begin{equation}
                                        \label{eq11.30.1}
u_{x^n}(0){\bar \psi}_{x^jx^j}(0)+u_{x^jx^j}(0)=0,
\end{equation}
which along with
Theorem \ref{theorem3.2} gives a bound for
$u_{x^jx^j}(0)$.

Next,
let's estimate the mixed derivatives $u_{x^jx^n}(0)$. We start
with
introduce a few more objects. Denote $\bA$ to be the space of
all
skew-symmetric matrices and for $p\in \bA$ we
set
$$
a(w,p)=a(e^pwe^{-p})=e^pa(w)e^{-p},\,\,\sigma=\sqrt{2a},\,\,f(w,p,x)=f(w,x).
$$
For $\xi\in \bR^d$, we also define $P(w,x)\xi=P(x)\xi$ with
value in
$\bA$ by the
formula
$$
[P\xi]_{ij}=\psi_{1x^i(\xi)}\psi_{1x^j}-\psi_{1x^i}\psi_{1x^j(\xi)}.
$$

Since
$e^pC_me^{-p}=C_m$ and $f(e^pwe^{-p},x)=f(w,x)$, we have
$$
0=\inf_{w\in C_m}\{a^{ij}(w)u_{x^ix^j}(x)+f(w,x)\}
$$
$$=\inf_{w\in
C_m, p\in \bA}\{a^{ij}(w,p)u_{x^ix^j}(x)+f(w,p,x)\}.
$$

Owing to the proof of  Theorem 5.9 of \cite{Kr95}, there exist
positive numbers $\delta_1$ and $\delta_2$ depending only on $K$ and
$d$ such that by taking $B_1:=\delta_1I_d$ the following assumption
is satisfied. This assumption is exactly Assumption 1.2 (d) of
\cite{athm} with $K$ there equal to $0$.
\begin{assumption}
                      \label{assumption6.0.1}
For any $x\in
\partial D$, $\xi\bot \psi_{1x}(x)$, $|\xi|=1$, $w\in C_m$, $p=0$ we
have $(B_1\xi,\xi)=\delta_1$ and
$$
L^w\psi_1+a^{ij}B_{1ij}\leq -\delta_2
$$
$$
(B_1\xi,\xi)L^w\psi+\sum_k(\partial(\xi)\psi_{1(\sigma^k)})^2
+2(B_1\xi,\sigma^k)\partial(\xi)\psi_{1(\sigma^k)}\leq -\delta_2,
$$
where
$$
\partial(\xi)\psi_{1(\sigma^k)}=\frac{d}{dh}\psi_{1x^i}(x+h\xi)
\sigma^{ik}(w,hP(w,x)\xi)|_{h=0}.
$$
\end{assumption}

The estimation of the mixed second order derivatives is a direct
application of Theorem 1.10 of \cite{athm}.  First
notice that, as we mentioned before, \eqref{eq2.2.b} is equivalent to
\eqref{eq05.02.11}, which is uniformly elliptic. Next, to estimate
the mixed second order derivatives we consider the function ${\bar
u}=u/\psi$, which satisfies a higher dimensional
  elliptic Bellman equation on an auxiliary manifold. Then
the problem is reduced to the estimation of tangential first order
derivatives of ${\bar u}$ on the manifold. In turn, actually it
suffices to have $f(w,\cdot)$ to be in $C^{1}({\bar D})$, which is
already satisfied in our case due to \eqref{eq2.0.0.4}. As a
conclusion we get

\begin{lemma}
Under our assumptions, there exist positive constants $\rho=\rho(K,d)$ and
$N=N(K,d)$ such that for any $x\in \partial D$ and unit $\tau\bot
\psi_{x}(x)$ we have
\begin{equation*}
                %          \label{eq6.0.2}
|u_{(\tau)(n)}(x)|\leq N(1+\max_{\partial D(\rho)}(|u|+|u_x|)].
\end{equation*}
\end{lemma}
This immediately implies the estimate of $u_{x^jx^n}(0)$,
$j=1,...,n-1$.

We use the equation \eqref{eq1.0.2} itself to estimate the normal
second order derivative $u_{x^nx^n}(0)$. Equation \eqref{eq1.0.2} at
the origin can be rewritten as
\begin{equation}
                          \label{eq8.1}
u_{x^nx^n}(0)P_{m-1,d-1}\big(u_{x^1x^1}(0),\cdots,u_{x^{n-1}x^{n-1}}(0)\big)+G=g^{m-1}(0),
\end{equation}
where
$G$ is a sum of products of $u_{x^jx^j}(0)$ and
$u_{x^jx^n}(0)$,
$j=1,\cdots,n-1$. By the results of Section
\ref{mixed}, we
have
\begin{equation}
                         \label{eq8.3}
|G|\leq
N(K,d).
\end{equation}
Due to \eqref{eqm-convex},
\eqref{eq11.30.1}
and Lemma \ref{lemma5.5}, we get
$$
P_{m-1,d-1}\big(u_{x^1x^1}(0),\cdots,u_{x^{n-1}x^{n-1}}(0)\big)
$$
$$
=(-u_{x^n}(0))^{m-1}P_{m-1,d-1}\big({\bar
\psi}_{x^1x^1}(0),\cdots,{\bar \psi}_{x^{n-1}x^{n-1}}(0)\big)
$$
\begin{equation}
                          \label{eq8.2}
\geq \delta(D,K,d)P_{m-1,d-1}(\kappa^1,\cdots,\kappa^{d-1}) \geq
\delta(D,K,d)>0.
\end{equation}
Combining \eqref{eq8.1}, \eqref{eq8.3} and \eqref{eq8.2} together
yields
\begin{equation}
                          \label{eq8.4}
|u_{x^nx^n}(0)|\leq N(D,K,d).
\end{equation}

Thus for the non-degenerate case $g>0$ in $D$, we get the estimate
for $|u_{xx}|$ on $\partial D$, and subsequently in $D$ by Theorem
\ref{theorem2.1}. For the general case, we only have to use Lemma
\ref{lemma02.02}.
\begin{remark}
                                  \label{remark8.1}
It is worth noting that when we estimate the second order
derivatives on the boundary $\partial D$, we only use the fact
$$
\|f(w,\cdot)\|_{C^1({\bar D})}\leq N(K,d).
$$ So to get the boundary estimates, it
suffices to assume $g^{2(m-1)/m}\in C^1({\bar D})$, which in general
is weaker than the condition that $g$ itself is in $C^1({\bar D})$.
\end{remark}

\begin{remark}
                                  \label{remark8.2}
Our method can be carried over to a larger class of Hessian
equations
\begin{equation}
                                  \label{eq05.1.9}
P_m(u_{xx}+h)=g^{m-1}
\end{equation}
with zero boundary condition, where $h\in C^{1,1}(\bar
D,\bar{C}_{m})$ is a $\bar{C}_{m}$-valued function satisfying the
following condition:
$$
\text{Tr}\, h\leq\frac{md}{d-m+1}\frac{P_{m}^{1-1/m}(I)}
{P_{m-1}(I)}g^{1-1/m},
$$
$$
\text{Tr}\, h\not\equiv\frac{md}{d-m+1}\frac{P_{m}^{1-1/m}(I)}
{P_{m-1}(I)}g^{1-1/m}.
$$
Naturally, we look for solutions such that $u_{xx}+h\in\bar{C}_{m}$
(a.e.).

Indeed, as we mentioned before, for large $t$ on $\partial D$ we
have $-\psi_{1,xx}\in C_m$. Due to Corollary \ref{cor4.2} on
$\partial D$ we have $\text{diag}\{\kappa^1,\cdots,\kappa^{d-1}\}\in
C_{m-1}$. As before, we can rewrite \eqref{eq05.1.9} as a Bellman
equation
\begin{equation}
                                          \label{eq05.1.9.1}
\inf_{w\in C_m}\{L^{w}u(x)+f(w,x)+\text{Tr}(a(w)h(x))\}=0,
\end{equation}
where $L^w$ and $f(w,x)$ are defined in the same way as in
\eqref{eq2.2.b}. By the same method we can get the estimates of $u$
and $u_{x}$, and reduce the interior estimate of $u_{xx}$ to the
estimates of $u_{xx}$ on the boundary.

  Under the linear transformation introduced at the
beginning of Section \ref{mixed}, $h(0)$ becomes $\bar{h}(0)$.
Denote $\bar{h}_n$ to be the $(d-1)\times (d-1)$ submatrix obtained
by deleting the $n$th row and $n$th column of $\bar{h}(0)$. By
Corollary \ref{cor4.2}, $\bar{h}_n$ is in ${\bar C}_{m-1,d-1}$.
After estimating the tangential and mixed second order derivatives
on the boundary in a similar way, we can obtain the estimate of
$u_{x^nx^n}(0)$ by using the equation \eqref{eq05.1.9} itself and
the inequality
$$
P_{m-1,d-1}\big(\text{diag}\{u_{x^1x^1}(0),
\cdots,u_{x^{n-1}x^{n-1}}(0)\}+{\bar h}_n\big)
$$
$$
\geq P_{m-1,d-1}\big(\text{diag}\{u_{x^1x^1}(0),
\cdots,u_{x^{n-1}x^{n-1}}(0)\}\big)
$$
$$
=P_{m-1,d-1}\big(\text{diag}\{-u_{x^n}(0){\bar
\psi}_{x^1x^1}(0),\cdots,-u_{x^n}(0){\bar
\psi}_{x^{n-1}x^{n-1}}(0)\}\big)
$$
$$
\geq (-u_{x^n}(0))^{m-1} P_{m-1,d-1}\big({\bar
\psi}_{x^1x^1}(0),\cdots,{\bar \psi}_{x^{n-1}x^{n-1}}(0)\big) \geq
\delta(D,K,d)>0.
$$
Here in the first step we use the inequality
$$
P_{m-1}(A+B)\geq P_{m-1}(A)
$$ for any $A,B\in {\bar C_{m-1}}$ (cf. Lemma \ref{lemma4.2.1} or \eqref{eq12.10}).
\end{remark}

\mysection{Proof of Lemma \ref{lemma02.02}}
                     \label{prooflemma}
Let $g_n$ be a sequence of strictly positive functions in
$C^2(\bar D)$ such that for $n=1,2,\cdots$, the functions
$g_n(x)+2K|x|^2$ are convex on $\bar D$ and
$$\inf_Dg_n\geq 1/(2n),\quad\|g_n\|_{C^1(\bar D)}\leq 2K,
\quad \|g_n-g\|_{C^0(\bar D)}\leq 1/n.$$ By our assumption there
exists $N=N(D,d,K)$ such that
\begin{equation}
                                 \label{eq02.02.3}
\|v_n\|_{C^2({\bar \Omega})}\leq N(D,d,K),\quad m=1,2,\cdots,
\end{equation}
where $v_n\in C^2({\bar D})\cap C^4(D)$ is the solution of
\eqref{dirichlet} with $g_n$ in place of $g$ and with zero Dirichlet
boundary condition. By Arzel\`a-Ascoli theorem, after passing to a
subsequence if necessary, $\{v_n\}$ converges in $C^1(\bar D)$ to a
function $v$. Again by \eqref{eq02.02.3} we get $v\in C^{1,1}(\bar
D)$ and
$$
\|v\|_{C^{1,1}({\bar \Omega})}\leq N(D,d,K).
$$
Owing to the uniqueness of the admissible weak solution, we obtain
$$u=v,\quad \|u\|_{C^{1,1}({\bar \Omega})}\leq N(D,d,K).$$ This
completes the proof of (i).

To prove (ii) we use the idea in the proof of Lemma 7.3.4 \cite{nonlinear}.
Let $\Omega$ be a countable dense subset of $C_m$. Obviously, one
has
$$
\inf_{w\in C_m}[L^wu+f]=\inf_{w\in \Omega}[L^wu+f].
$$

For any $w\in \Omega$, in $D$ it holds that
\begin{equation}
                                     \label{eq02.05.1}
L^wv_n(x)+f_n(w,x)\geq 0,
\end{equation}
where
$$
f_n(w,x)=-m(d-m+1)^{-1}P_m^{1-1/m}(w)P_{m-1}^{-1}(w)g_n^{1-1/m}(x).
$$
After multiplying \eqref{eq02.05.1} by a nonnegative function
$\eta\in C^{\infty}_0(D)$, integrating by parts, passing to the
limit over the sequence $n$ and integrating by parts again, we
obtain
$$
\int_D[a^{ij}(w)\eta(x)u_{x^ix^j}(x)+\eta(x)f(w,x)]\,dx\geq 0.
$$
Because $\eta\in C^{\infty}_0(D)$ is arbitrary, we further get
$$
a^{ij}(w)u_{x^ix^j}(\cdot)+f(w,\cdot)\geq 0\quad \text{a.e. in}\,D.
$$
Since $\Omega$ is a countable set, we reach
\begin{equation}
                                     \label{eq02.05.02}
\sup_{w\in \Omega}[Lu(x)+f(w,x)]\geq 0\quad \text{a.e. in}\,D.
\end{equation}

Next we prove the opposite inequality. Here we use again the method
by which Lemma 7.3.4 of \cite{nonlinear} is proved. Let $\epsilon<1$
be a positive number. Recall that $v_{n}$ satisfies the Bellman
equation
$$
\inf_{w\in C_m}\{L^{w}v_n(x)+f_n(w,x)\}=0.
$$
Given any $n$ such that $1/n<\epsilon$, we have
$$
\inf_{w\in C_m}\{L^{w}v_n(x)+\epsilon \Delta v_n(x)+f(w,x)\}
$$
$$
\leq \sup_{w\in C_m}\{\epsilon \Delta v_n(x)+f(w,x)-f_n(w,x)\} \leq
N(D,d,K)\epsilon,
$$
$$
\inf_{w\in C_m}\{L^{w}v_n(x)+\epsilon \Delta
v_n(x)+f(w,x)-N(D,d,K)\epsilon\}\leq 0.
$$
Note that the elliptic operator $\inf_{w\in C_m}[L^{w}+\epsilon
\Delta]$ is uniformly non-degenerate. Owning to Theorem 3.6.3
\cite{nonlinear}, after passing to the limit over the sequence $n$
we obtain
$$
\inf_{w\in C_m}\{L^{w}u(x)+\epsilon \Delta
u(x)+f(w,x)-N(D,d,K)\epsilon\}\leq 0 \,\,\text{a.e. in}\,D.
$$
Letting $\epsilon\downarrow 0$ yields an inequality opposite to
\eqref{eq02.05.02}. The lemma is proved.

\end{document}